\documentclass{elsart3-1}

\usepackage{amssymb}
\usepackage{amsmath,multicol}
\usepackage[english,francais]{babel}

\newtheorem{theorem}{Theorem}[section]

\newtheorem{e-proposition}[theorem]{Proposition}
\newtheorem{corollary}[theorem]{Corollary}
\newtheorem{e-definition}[theorem]{Definition}
\newtheorem{e-remark}[theorem]{Remark}

\renewcommand{\theequation}{\arabic{equation}}
\def\C{\bb C} 
\def\N{\bb N}
\def\Z{\bb Z}

\def\go{\mathfrak}
\def\bb{\mathbb}
\def\cal{\mathcal}

\def\og{\leavevmode\raise.3ex\hbox{$\scriptscriptstyle\langle\!\langle$~}}
\def\fg{\leavevmode\raise.3ex\hbox{~$\!\scriptscriptstyle\,\rangle\!\rangle$}}

\journal{the Acad\'emie des sciences}
\begin{document}

\centerline{Representation Theory}
\begin{frontmatter}


\selectlanguage{english}
\title{Algebras of invariant differential operators on a class of  multiplicity free spaces.}


\selectlanguage{english}
\author[authorlabel1]{Hubert Rubenthaler},
\ead{rubenth@math.u-strasbg.fr}

\address[authorlabel1]{Institut de Recherche Math\'ematique Avanc\'ee, Universit\'e de Strasbourg et CNRS, 7 rue Ren\'e Descartes, 67084 Strasbourg cedex}


\medskip
\begin{center}
{\small Received *****; accepted after revision +++++\\
Presented by £££££}
\end{center}
\begin{abstract}
\selectlanguage{english}
Let $G$ be a connected reductive algebraic group and let $G'=[G,G]$ be its derived subgroup. Let $(G,V)$ be a  multiplicity free representation with a one dimensional quotient (see definition below). We prove that the algebra $D(V)^{G'}$ of $G'$-invariant differential operators with polynomial coefficients on V, is a quotient of   a so-called  Smith  algebra over its center.  Over ${\bb C}$ this class of algebras was introduced by S.P. Smith \cite{Smith} as a class of algebras similar to ${\cal U}({\go s}{\go l}_{2})$.    Our result generalizes the case of the Weil representation,     where the associative algebra generated by $Q(x)$ and $Q(\partial)$ ($Q$ being a non degenerate quadratic form on $V$) is a quotient of ${\cal U}({\go s}{\go l}_{2})$. Other structure results are obtained when $(G,V)$ is a regular prehomogeneous vector space of commutative parabolic type.
{\it To cite this article: Hubert Rubenthaler,   C. R. Acad. Sci. Paris, Ser. I 340 (2008).}

\vskip 0.5\baselineskip

\selectlanguage{francais}
\noindent{\bf R\'esum\'e} \vskip 0.5\baselineskip \noindent
{\bf Alg\`ebres d'op\'erateurs diff\'erentiels invariants associ\'es \`a certains espaces sans multiplicit\'es. }
Soit $G$ un groupe alg\'ebrique r\'eductif connexe et soit $G'=[G,G]$ son groupe d\'eriv\'e. Soit $(G,V)$ un espace sans multiplicit\'es ayant un quotient unidimensionel (voir la d\'efinition ci-dessous). Nous montrons que l'alg\`ebre  $D(V)^{G'}$ des op\'erateurs diff\'erentiels  \`a coefficients poynomiaux $G'$-invariants sur $V$, est isomorphe \`a un quotient d'une alg\`ebre de Smith sur son centre. Sur ${\bb C}$  cette classe d'alg\`ebres,   avait  \'et\'e introduite par S.P.  Smith (\cite{Smith}) comme une classe d'alg\`ebres semblables  \`a  ${\cal U}({\go s}{\go l}_{2})$. Notre r\'esultat g\'en\'eralise le cas de la repr\'esentation de Weil, o\`u l'alg\`ebre associative engendr\'ee par $Q(x)$ et $Q(\partial)$ ($Q$ \'etant une forme quadratique non d\'eg\'en\'er\'ee sur $V)$, est un quotient de ${\cal U}({\go s}{\go l}_{2})$. D'autres r\'esultats de structure   sont obtenus lorsque $(G,V)$ est un espace pr\'ehomog\`ene parabolique commutatif r\'egulier.
{\it Pour citer cet article~: Hubert Rubenthaler,  C. R. Acad. Sci.
Paris, Ser. I 340 (2008).}

\end{abstract}
\end{frontmatter}
 \selectlanguage{english}
\section{Introduction} 
For any smooth affine algebraic variety $M$, we shall denote by $\C[M]$ the algebra of regular functions on $M$, and by ${\bf D}(M)$ the algebra of differential operators on $M$. Moreover, if an algebraic group $H$ acts on $M$, we will denote respectively by $\C[M]^H$ and ${\bf D}(M)^H$ the algebras of $H$-invariant  regular functions and of $H$-invariant differential operators on $M$.

Let $Q(x)$ be a non degenerate quadratic form on ${\bb C}^n$ and let $Q(\partial)$ be the differential operator with constant coefficients defined by the condition $(Q(\partial)Q)(0)=1$. It is well known that these two differential operators generate a Lie subalgebra of ${\bf D}({\bb C}^n)$ which is isomorphic to ${\go s}{\go l}_{2}$ (note that this Lie algebra completely determines the infinitesimal Weil representation of ${\go {sp}}(n,\C)$). Consequently the associative subalgebra ${\bb C}[Q(x),Q(\partial)]$ of ${\bf D}(\C^n)$, generated by these two operators, is a quotient of ${\cal U}({\go s}{\go l}_{2})$. A natural question one can ask   is the following: what happens if $Q(x)$ is replaced by an arbitrary homogeneous polynomial $P$ of degree $n$? J.I. Igusa (\cite{Igusa}) has proved that if $\partial ^{\circ}P\geq 3$, then the Lie algebra generated by $P(x)$ and $P(\partial)$ is infinite dimensional.  In this paper we basically study the case where $P(x)=\Delta_{0}(x)$ is the fundamental relative invariant of a so-called multiplicity-free space $(G,V)$ with one dimensional quotient (\cite{Levasseur}, see \ref{def-multiplicityfree-1-quotient} below).  Among several results concerning various algebras of differential operators, our main result is Theorem \ref{th-principal} below which says that the algebra ${\bf D}(V)^{G'}$ is a quotient of a so-called Smith algebra over its center. As Smith algebras are "similar" to ${\cal U}({\go {sl}}_{2})$, this can be viewed as a generalisation of the infinitesimal Weil representation. The stucture of the algebra of radial components of ${\bf D}(V)^{G'}$ is obtained as a corollary.
\section{Invariant differential operators on multiplicity free spaces with one dimensional quotient.} Let $G$ be an  reductive algebraic   group over ${\bb C}$, and let $G'=[G,G]$ be its derived group. Let $(G,V)$ be a rational finite dimensional linear representation of $G$. Recall that $(G,V)$ is said to be {\it multiplicity free} if the associated  representation of $G$ on ${\bb C}[V]$ decomposes without multiplicities. For the classification and properties   of multiplicity free spaces we refer to \cite{Benson-Ratcliff}, \cite{Leahy}, \cite{Knop}. Let us first recall that a multiplicity free space  is necessarily a prehomogeneous vector space (abbreviated $PV$) not only under $G$ but also under any Borel subgroup $B$ of $G$.
\begin{e-definition}\label{def-multiplicityfree-1-quotient}{\rm (Levasseur \cite{Levasseur})}
 A multiplicity-free space $(G,V)$ is said to have a one dimensional quotient if there exists  a polynomial $\Delta_{0}$ such that $\Delta_0 \notin \C[V]^G$ and such that $\C[V]^{G'}=\C[\Delta_{0}]$. 
  \end{e-definition} 
 From now on we assume that $(G,V)$ is a multiplicity free space with a one dimensional quotient. We remark first that the polynomial $\Delta_{0}$ is necessarily irreducible. There exists then a non trivial character $\chi_{0}$ of $G$, such that $\forall g\in G, \forall x\in V, \, \Delta_{0}(gx)=\chi_{0}(g)\Delta_{0}(x)$. In other words $\Delta_{0}$ is a relative invariant of the $PV$ ($G,V$), with character $\chi_{0}$.   We define $\Omega=\{x\in V\,|\, \Delta_{0}(x)\neq0\}$. The algebra of regular functions on $\Omega$ is given by $\displaystyle \C[\Omega]= \{f=\frac{P}{\Delta_{0}^k}\,|\, P\in \C[V], k\in \N\}$. From the definitions it is easy to see that the natural action of $G$ on $\C[\Omega]$ is still without multiplicities. As an easy consequence of this fact, one can prove the following result.
 
 \begin{e-proposition}\label{prop-coef-poly}
 One has ${\bf D}(\Omega)^G= {\bf D}(V)^G$ (that is, any $G$-invariant differential operator with regular coefficients on $\Omega$ has polynomial coefficients). Moreover   this algebra is commutative.
 
 \end{e-proposition}
 Hence we have the following inclusions among the algebras ${\bf D}(V)^G  = {\bf D}(\Omega)^G , {\bf D}(V)^{G'},    {\bf D}(\Omega)^{G'}$:
  $$\begin{matrix}{\bf D}(V)^G & = &{\bf D}(\Omega)^G\\
\downarrow & { }&\downarrow \\
{\bf D}(V)^{G'} & \longrightarrow &  {\bf D}(\Omega)^{G'}.
\end{matrix}$$
In the case where $(G,V)$ is a regular $PV$ of commutative parabolic type, or equivalently when $V$ is  a simple Jordan algebra over $\C$ and $G$ its structure group,  the  result of Proposition \ref{prop-coef-poly} was essentially known by computing explicit generators of ${\bf D}(\Omega)^G$ (\cite{Nomura}, \cite{Yan}).
 
 Let us define:\hskip 10pt ${\cal T}_{0}={\bf D}(V)^G= {\bf D}(\Omega)^G,\quad {\cal T}= {\bf D}(\Omega)^{G'}. $  
 
 Let us also set:
  $$X=\Delta_{0}(x), \qquad X^{-1}= \frac{1}{\Delta_{0}(x)}, \qquad Y=\Delta_{0}(\partial)\qquad E= {\text {\rm  Euler operator}}$$
 where $\Delta_{0}(\partial)$ is the differential operator with constant coefficients on $V$ defined by the condition $(\Delta_{0}(\partial)\Delta_{0})(0)=1$. Clearly $E\in {\cal T}_{0}={\bf D}(V)^G$, $X,Y\in {\bf D}(V)^{G'}$, and $X^{-1}\in {\bf D}(\Omega)^{G'}={\cal T}$.

 \begin{e-definition}\label{def-tau}
 The automorphism $\tau$ of the algebra ${\cal T}$ is defined by
 $$\forall D\in {\cal T},\qquad \tau(D)=XDX^{-1} .$$
 \end{e-definition}

 It is easy to see that ${\cal T}_{0}$ is stable under $\tau$. The preceding definition can also be read $XD=\tau(D)X$ for $D\in {\cal T}$. In fact a similar relation  is true for $Y$: for all $D\in {\cal T}$, we have $DY=Y\tau(D)$.

 Let us denote by ${\cal T}_{0}[X,Y]$ the subalgebra of ${\cal T}$ generated by ${\cal T}_{0}$, $X$, and $Y$ and by ${\cal T}_{0}[X,X^{-1}]$ the subalgebra of ${\cal T}$ generated by ${\cal T}_{0}$, $X$, and $X^{-1}$.

 \begin{e-proposition} \label{prop-graduation}Let $d_{0}$ be the degree of $\Delta_{0}. $We have:
 
\noindent {\rm 1)} $\displaystyle {\bf D}(V)^{G'}={\cal T}_{0}[X,Y]= \oplus_{p\in \N}{\cal T}_{0}X^p \bigoplus \oplus_{k\in \N^*}{\cal T}_{0}Y^k$, and moreover 
 ${\cal T}_{0}X^p={\cal T}_{0}[X,Y]_{p} $
 
 \noindent$=\{D\in {\cal T}_{0}[X,Y]\,|\, [E,D]=pd_{0}D\},  p\in \N$, ${\cal T}_{0}Y^k={\cal T}_{0}[X,Y]_{k}=\{D\in {\cal T}_{0}[X,Y]\,|\, [E,D]=-kd_{0}D\}, k\in \N^*.$ 
 
\noindent   {\rm 2)} ${\cal T}=  {\bf D}(\Omega)^{G'}={\cal T}_{0}[X,X^{-1}]=\oplus _{p\in \Z}{\cal T}_{0}X^p$, and moreover 
  ${\cal T}_{0}X^p={\cal T}_{p}=\{D\in {\cal T}\,|\, [E,D]=pd_{0}D \}, p\in \Z$.
      \end{e-proposition}

      \begin{e-remark}\label{remark-du-a-levasseur}
      Assertion  {\rm 1)} is due to Levasseur {\rm(\cite{Levasseur}, Th. 4.11.)}. In the particular case of $PV$'s of commutative parabolic type the result was obtained independently by the author {\rm(\cite{Rubenthaler-arxiv})}.
    \end{e-remark}
     
     Let now ${\cal Z}({\cal T})$ be the center of ${\cal T}$. It can be shown that ${\cal Z}({\cal T})$ is also the center of ${\bf D}(V)^{G'}={\cal T}_{0}[X,Y]$ and that  $D\in {\cal Z}({\cal T})\Longleftrightarrow D \in {\cal T}_{0} \text{  and  }  \tau(D)=D$).
    
  The key result is then  the following:
  \vskip 3pt
  \begin{theorem}\label{th-structure-T0}\hfill
 
\noindent   {\rm 1)}  ${\cal T}_{0}={\cal Z}({\cal T})\oplus E{\cal T}_{0}$
  
 \noindent  {\rm 2)}  Any element $H\in {\cal T}_{0}$ can be uniquely  written  in the form  

\centerline {$ H=H_{0}+EH_{1}+\dots+E^kH_{k}, \, H_{i}\in {\cal Z}({\cal T}),\, k\in \N$ }

  \end{theorem}
{Sketch of proof:}
   Let $B$ be a Borel subgroup of $G$.   As mentioned before the space $(B,V)$ is prehomogeneous. Let $\Delta_{0},\Delta_{1}, \dots, \Delta_{r}$ be the set of fundamental relative invariants of $(B,V)$, where $\Delta_{0}$ is the unique fundamental relative invariant of $(G,V)$ (the unicity is due to the one dimensional quotient hypothesis we make here). We denote by $d_{i}$ (resp. $\lambda_{i}$) the degree (resp. the infinitesimal character) of $\Delta_{i}$ $(i=0,\dots,r$). Let ${\go b}$ be the Lie algebra of $B$. Let ${\go t}\subset {\go b}$ be a Cartan sugalgebra of ${\go g}$, and let $\Sigma$ be the set of roots of $({\go g}, {\go t})$. Denote by $W$ the Weyl group of $\Sigma$.  Let $\Sigma^+$ be the set of positive roots such that ${\go b}={\go t}\oplus \Sigma_{\alpha\in \Sigma^+}{\go g}^{\alpha}$. Let ${\go a}^*=\oplus_{i=0}^r \C \lambda_{i}\subset {\go t}^*$ and let $\rho=\frac{1}{2}\Sigma_{\alpha\in \Sigma^+}\alpha$. Define $A={\go a}^*+\rho\subset{\go t}^*$. F. Knop (\cite {Knop}, \cite {Benson-Ratcliff}) has proved that there exists an isomorphism (the so-called Harish Chandra isomorphism) $h:{\cal T}_{0}={\bf D}(V)^{G}\longrightarrow \C[A]^{W_{0}}$, where $W_{0}$ is a subgroup of the stabilizer of $A$ in $W$. If $\lambda=\Sigma_{i=0}^{r}\alpha_{i}\lambda_{i}\in {\go a}^*$, then $h(E)(\lambda+\rho)=\Sigma_{i=0}^r\alpha_{i}d_{i}$. Let $M=\{\mu\in A\,|\, h(E)(\mu)=0\}$ and let $I(M)=\{P\in \C[A]^{W_{0}}\, | \,P_{|_{M}}=0\}$. We still denote by $\tau $ the automorphism of $\C[A]^{W_{0}}$ which corresponds to $\tau_{|_{{\cal T}_{0}}}$ via the Harish Chandra isomorphism $h$. Let $\C[A]^{W_{0}, \tau}$ be the set of $\tau$-invariants in $\C[A]^{W_{0}}$. We prove in fact that $\C[A]^{W_{0}}=\C[A]^{W_{0},\tau}\oplus I(M)$ and the pullback by $h $ of this decomposition proves the first assertion. The second assertion is a direct consequence of the first.
  
  \begin{e-definition}\label{def-algebre-smith}Let ${\bf A}$ be a commutative  associative algebra over $\C$, with unit element $1$ and without zero divisors. Let $f, u \in {\bf A}[t]$ be two polynomials  in one variable with coefficients in ${\bf A}$. Let $n\in \N^*$.
  
\noindent  {\rm 1)} {\rm {(Smith \cite{Smith}}, over ${\bb C}$)}The (generalized) Smith algebra $S({\bf A}, f, n)$ is the associative algebra over ${\bf A}$ with generators $(x,y, e)$ subject to the relations 
  
  \centerline{$[e,x]= nx,\quad [e,y]=-ny, \quad [y,x]=f(e)$}
  
  \noindent  {\rm 2)} The algebra $U({\bf A}, u, n)$ is the associative algebra over ${\bf A}$ with generators $\tilde x, \tilde y, \tilde e$ subject to the relations 
  
  \centerline{$[\tilde e,\tilde x]=n\tilde x,\quad [\tilde e, \tilde y]=-n\tilde y, \quad \tilde x\tilde y=u(\tilde e), \quad \tilde y\tilde x=u(\tilde e + n)$}
  \end{e-definition}
  
 \noindent  One can show  that if $f\in {\bf A}[t]$ is given, and if $u\in {\bf A}[t]$ is the unique polynomial (up to scalars in ${\bf A}$) given by $f(t)=u(t+n)-u(t)$, then $U({\bf A}, u, n)\simeq S({\bf A}, f,n)/(\Omega)$, where $\Omega= xy-u(e)$ generates the center of $S({\bf A}, f, n)$ over ${\bf A}$.
 
As $XY\in {\cal T}_{0}$, we know from Theorem \ref{th-structure-T0} that there exists a unique $u_{XY}\in {\cal Z}({\cal T})[t]$ such that $XY=u_{XY}(E)$. Our main result is then the following:
\begin{theorem}\label{th-principal}The mapping {$\tilde x \longmapsto X, \quad \tilde y \longmapsto Y, \quad \tilde e \longmapsto E$} extends uniquely to an isomorphism of ${\cal Z}({\cal T})$-algebras between $U({\cal Z}({\cal T}), u_{XY}, d_{0})$ and ${\cal T}_{0}[X,Y]={\bf D}(V)^{G'}$. 
\end{theorem}

{Sketch of proof:} From Theorem \ref{th-structure-T0} we know that ${\cal T}_{0}[X,Y]={\bf D}(V)^{G'}$ is generated over ${\cal Z}({\cal T})$ by $X,Y,E$. As these three elements satisfy the defining relations of $U({\cal Z}({\cal T}), u_{XY}, d_{0})$, there exists a surjective morphism of ${\cal Z}({\cal T})$-algebras $\varphi: U({\cal Z}({\cal T}), u_{XY}, d_{0})\longrightarrow {\cal T}_{0}[X,Y]$. Theorem \ref{th-structure-T0} is used again to prove that any $D\in {\cal T}_{0}[X,Y]$ can be written uniquely in the form $D=\sum_{\ell>0,k\geq 0}\alpha_{k,\ell}Y^{\ell}E^k+\sum_{m\geq0,r\geq0}\beta_{m,r}X^mE^r$ with $\alpha_{k,\ell},\beta_{m,r}\in {\cal Z}({\cal T})$. This implies that $\varphi$ is bijective.

\noindent    We define the radial component $\overline D$ of an operator  $D\in {\cal T}_{0}[X,Y]={\bf D}(V)^{G'}$, by $D(\psi\circ \Delta_{0})= (\overline D \psi)\circ \Delta_{0}$, for $\psi \in \C[t]$. If $D\in  {\cal Z}({\cal T}) $, then $\overline D$ is just a constant. If $u(t)=D_{0}+D_{1}t+\dots+D_{k}t^k\in {\cal Z}({\cal T})[t]$, we define $\overline u \in \C[t]$, by $\overline u(t)=\overline D_{0}+\overline D_{1}t+\dots+\overline D_{k}t^k $. Let us also denote by $R$ the image of the map $D\longmapsto \overline D$ (the so-called algebra of radial components).  As a consequence we obtain a new proof of the following result. 
\begin{corollary}{\rm (Levasseur \cite{Levasseur})} $R\simeq U(\C, \overline u_{XY}, d_{0})$.
\end{corollary} 

\section{The case of irreducible regular $PV$'s of parabolic commutative case}
Let $\tilde{\go g}= V^-\oplus {\go g}\oplus V^+$ be a 3-grading of a simple Lie algebra $\tilde{\go g}$ over $\C$ ($[{\go g}, V^\pm]\subset V^\pm, [V^-,V^+]\subset {\go g}$). Then if  $G$ is the connected subgroup of the adjoint group of $\tilde {\go g} $ with Lie algebra ${\go g}$, it is well known that $(G,V^+)$ is a multiplicity free space (these $PV$'s  are said to be of {\it commutative parabolic type}, they are in one to one correspondence with Jordan triple systems or with hermitian symmetric spaces). We assume further that  $(G,V^+)$ is regular as a $PV$ (these objects are in one to one correspondence with simple Jordan algebras, or with hermitian symmetric spaces of tube type). These regular $PV$'s of commutative parabolic type form a subclass of the multiplicity free spaces with one dimensional quotient. For these we can give some further information concerning the preceding algebras of differential operators. First of all, the number of $\Delta_{i}$'s is equal to $d_{0}$. In the notations of  section 2, $\Delta_{0}, \Delta_{1}, \dots, \Delta_{d_{0}-1}$ is the set of fundamental relative invariants of $(B,V^+)$. For simplicity define $n=d_{0}-1$. For ${\bf a}=(a_{0},\dots,a_{n})\in \N^{n+1}$, let $V_{\bf a}$ be the irreducible $G$-submodule of $\C[V^+]$ generated by $\Delta^{\bf a}=\Delta_{0}^{a_{0}}\dots\Delta_{n}^{a_{n}}$. It is well known that the $G$-module $V_{\bf a}\otimes V_{\bf a}^*$ contains up to constants, a unique $G$-invariant vector $C_{\bf a}$, which can be viewed as a $G$-invariant differential operator (see \cite{Howe-Umeda}). These operators are sometimes called Capelli operators. Define ${\bf 1}_{j}=(0,\dots,0,1,0,\dots,0)\in \N^{n+1}$. Note that $XY=C_{{\bf 1}_{0}}$ and that $E=C_{{\bf 1}_{n}}$.

\begin{e-proposition}\hfill

\noindent {\rm 1)} Let $\C[E, X, X^{-1}, Y] $ be the associative subalgebra of ${\bf D}(V^+)$ generated by $E, X, X^{-1}, Y$. We have $\C[E, X, X^{-1}, Y] ={\cal T}$.

\noindent {\rm 2)} { \rm (Yan, \cite{Yan})} We have ${\cal T}_{0}=\C[XY,C_{{\bf 1}_{1}},\dots, C_{{\bf 1}_{n}}]$,and $XY,C_{{\bf 1}_{1}},\dots, C_{{\bf 1}_{n}}$ are algebraically independent.

\noindent {\rm 3)} We have ${\cal T}_{0}[X,Y]={\bf D}(V^+)^{G'}=\C[X,Y,C_{{\bf 1}_{1}},\dots, C_{{\bf 1}_{n}}]$.
\end{e-proposition}






\end{document}